\documentclass[11pt]{article}
\usepackage{enumerate}
\usepackage{amssymb,a4wide,latexsym,makeidx,epsfig,fleqn}
\usepackage{amsthm}
\usepackage{amsmath}
\usepackage{enumerate}

\newtheorem{theorem}{Theorem}[section]

\newtheorem{lemma}[theorem]{Lemma}

\begin{document}
\textwidth 150mm \textheight 225mm
\title{The rank of a signed graph in terms of girth \footnote{This work is supported by the National Natural Science Foundations of China (No. 11901253), the Natural Science Foundation for Colleges and Universities in Jiangsu Province of China (No. 19KJB110009), and the Science Foundation of Jiangsu Normal University (
No. 18XLRX021).}}
\author{{Yong Lu, Qi Wu\footnote{Corresponding author.}}\\
{\small  School of Mathematics and Statistics, Jiangsu Normal University,}\\ {\small  Xuzhou, Jiangsu 221116,
People's Republic
of China.}\\
{\small E-mail: luyong@jsnu.edu.com, wuqimath@163.com}}

\date{}
\maketitle
\begin{center}
\begin{minipage}{120mm}
\vskip 0.3cm
\begin{center}
{\small {\bf Abstract}}
\end{center}
{\small Let $\Gamma=(G,\sigma)$ be a signed graph and $A(G,\sigma)$ be its adjacency matrix. Denote by $gr(G)$ the girth of $G$, which is the length of the shortest cycle in $G$. Let $r(G,\sigma)$ be the rank of $(G,\sigma)$.
In this paper, we will prove that $r(G,\sigma)\geq gr(G)-2$ for a signed graph $(G,\sigma)$. Moreover, we characterize all extremal graphs which satisfy the equalities $r(G,\sigma)=gr(G)-2$ and $r(G,\sigma)=gr(G)$.

 }
\end{minipage}
\end{center}

\section{Introduction }



In this paper, we only consider the graphs which have no multiedges and loops. Let $G=(V(G),E(G))$ be a simple graph, where $V(G)$ and $E(G)$ are the vertex
set and the edge set of $G$, respectively. Let $V(G)=\{v_{1},v_{2},\ldots,v_{n}\}$,  the \emph{adjacency matrix} $A(G)$ of $G$ is the symmetric $n\times n$ matrix with entries $a_{ij}=1$ if $v_{i}v_{j}\in E(G)$ and $a_{ij}=0$ otherwise. Denote by $r(G)$(resp. $\eta(G)$) the \emph{rank}(resp. \emph{nullity}) of $G$, which is the multiplicity of nonzero (resp. 0) eigenvalue of $A(G)$. We use $N_{G}(x)$ to denote the \emph{neighbor set} of a vertex $x$ in $G$. Denote by $d_{G}(x)$ the \emph{degree} of $x$ in $G$ which equals $|N_{G}(x)|$. A \emph{pendant vertex} of graph $G$ is a vertex with degree one in $G$. Denote by $P_{n}$ and $C_{n}$ a path and cycle of order $n$, respectively. We call a subset $I$ of $V(G)$ an \emph{independent set} if $uv\notin E(G)$ for any two vertices $u,v$ of $I$
in a graph $G$. We use $K_{m,n}$ to denote the complete bipartite graph.

Let $S\subseteq V(G)$ and $S\neq \phi$, we use $G[S]$ to denote the \emph{induced subgraph} of $G$, whose vertex set is $S$ and edge set is the set of those edges of $G$ that have both ends in $S$.
 Denote by $G-S$ the induced subgraph obtained from $G$ by deleting each vertex in $S$ and its incident edges. For convenience, if $S=\{x\}$, we write $G-x$ instead of $G-\{x\}$. For the induced subgraph $H$ of $G$, we use $H+x$ to denote the subgraph of $G$ induced by the vertex set $V(H)\cup \{x\}$.

The \emph{girth} of $G$, denoted by $gr(G)$ (or simplify $g$), is the size of a shortest cycle of $G$. The \emph{length} of a path is the number of edges in the path. Denote by $d(x,y)$ the \emph{distance} between $x$ and $y$ which is the length of a shortest path between vertices $x,y$ of graph $G$. If $H$ is a subgraph of graph $G$ and $x$ is a vertex outside $H$ in $G$, denote by $d(x,H)$ the length of shortest path from $x$ to all vertices of $H$.

 A \emph{signed graph} $\Gamma=(G, \sigma)$ consists of a simple graph $G$ with edge set $E$ and a mapping $\sigma: E\rightarrow\{+, -\}$, where $G$ is called the \emph{underlying graph} of $(G, \sigma)$. The \emph{adjacency matrix}  of $(G, \sigma)$, denoted by $A(G, \sigma)=a^{\sigma}_{ij}=\sigma(v_{i}v_{j})a_{ij}$, where $a_{ij}\in A(G)$. We use $r(G, \sigma)$(resp. $\eta(G,\sigma)$) to denote the \emph{rank}(resp. \emph{nullity}) of a signed graph $(G, \sigma)$ which is \emph{rank}(resp. \emph{nullity}) of $A(G,\sigma)$.  For a signed graph $(G,\sigma)$ of order $n$, $r(G,\sigma)+\eta(G,\sigma)=n$. Denote by $N_{+}(v)$
(resp. $N_{-}(v)$) the set of vertices joining $v$ with positive
edges (resp. negative edges) in a signed graph $(G,\sigma)$. A signed graph $(G,\sigma)$ is called \emph{nonsingular} if the determinant of $A(G,\sigma)$ is nonzero.

 Let $(C_{n},\sigma)$ be a cycle of $(G,\sigma)$. The \emph{sign} of $(C_{n},\sigma)$, denoted by $sgn(C_{n},\sigma)$, is defined as $\prod_{e\in E(C_{n},\sigma)}\sigma(e)$. If $sgn(C_{n},\sigma)=+$ (or $sgn(C_{n},\sigma)=-$), then we say that $(C_{n},\sigma)$ is \emph{positive} (or \emph{negative}). If all the cycles  of $(G,\sigma)$ are positive, then $(G,\sigma)$ is \emph{balanced}, and \emph{unbalanced} otherwise.

Denote by $S_{n}$ the \emph{$n$-vertex star} ($n\geq 2$) which is the complete bipartite graph $K_{1,n-1}$. For $n\geq 3$, We call the vertex of degree $n-1$ in $S_{n}$ the \emph{star center} (For the star $S_{2}$, each
of its vertices is referred to as its star center).  A subgraph $G_{1}$ of a graph $G$ is
called a \emph{pendant star} of $G$ if $G_{1}$ is a star such that its star center is its only vertex that has a
neighbor not in $G_{1}$.

Let $c(G)=|E(G)|-|V (G)|+\omega(G)$ be the dimension of cycle space of a graph $G$, where $\omega(G)$ is the number of connected components of $G$. A connected graph $G$ is called \emph{unicyclic}, \emph{bicyclic} and \emph{tricyclic} if $c(G)=1$, $c(G)=2$ and $c(G)=3$, respectively. A \emph{canonical unicyclic graph} is a cycle with pedant stars attached at none, some or all of its vertices (see Figure 1).  Let $P_{p}, P_{l}, P_{q}$ be three paths, where $\min\{p, l, q\}\geq 2$ and at
most one of $p, l, q$ is $2$. Let $\theta(p, l, q)$  be the graph obtained from $P_{p}, P_{l}$
and $P_{q}$ by identifying the three initial vertices and terminal vertices.

Let $x$ and $y$ be two different vertices of a signed graph $(G,\sigma)$. If $N_{G}(x)=N_{G}(y)=V_{1}$
and there is a nonzero constant $k\in\{-1,+1\}$ such that $\frac{\sigma(xz)}{\sigma(yz)}=k$ for any vertex $z\in V_{1}$
, then $x$ is said to be a \emph{multiple} of $y$ in $(G,\sigma)$. Clearly, $r(G,\sigma)=r((G,\sigma)-x)=r((G,\sigma)-y)$, if $x$ is a multiple of $y$ in $(G,\sigma)$. Let $(H,\sigma)$ be the signed graph obtained from $(G,\sigma)$ by deleting multiples successively. If $(H,\sigma)$ has no two vertices satisfying that one is the multiple of the other, we call $(H,\sigma)$
 the \emph{reduced graph} of $(G,\sigma)$.

Collatz et al. \cite{CS} had wanted to obtain all graphs of order $n$ with $r(G)<n$. Until today, this problem is also unsolved.
In mathematics, the rank (or nullity) of a graph has attracted a lot of researchers' attention, they focus on the  bounds for  the rank (or nullity) of a simple graph $G$ \cite{FHLL, LGUO, MFANG, MWTDAM, RCZ, LWANG, WLMA, WW, ZWS, ZWT}.

In recent years, the researches on the rank and nullity of signed graphs obtained increased
attention. In \cite{LH}, Liu and You researched the nullity of signed graphs. Fan et al. studied the nullity of unicyclic signed graphs and bicyclic signed graphs in \cite{FWW} and \cite{FDD}, respectively. Belardo et al. investigated the spectral characterizations and the Laplacian coefficients of signed graphs in \cite{BP} and \cite{BS}, respectively. He et al. studied the relation between the rank of a signed graph
and the matching number of its underlying graph in \cite{HHL}. Li et al. \cite{LXWX} studied the rank of a signed graph
 in terms of independence number. Lu and Wu \cite{LUWU1} proved that no signed graph with the nullity $\eta(G,\sigma)=|V(G)|-2m(G)+2c(G)-1$, where $m(G)$ is the matching number of $G$. There are other researches of
a signed graph, readers can refer to those in \cite{LW, LWZ, WS, YFQ}.

 In this paper, we will study the relation among the rank of a signed graph $(G,\sigma)$ and
the girth of its underlying graph. For any connected signed graph $(G,\sigma)$, we will prove that $r(G,\sigma)\geq gr(G)-2$. Moreover, we characterize all extremal graphs which satisfy the equalities $r(G,\sigma)=gr(G)-2$ and $r(G,\sigma)=gr(G)$.  Our
results generalize the corresponding results on graph \cite{ZWT}.

 In Section 2, we give some lemmas about signed graphs. In Section 3, we characterize the relations between the rank of a signed graph and its girth.

\section{Preliminaries}
\quad For a signed graph $(G,\sigma)$, there has some lemmas.
\noindent\begin{lemma}\label{le:2.1}\cite{LWZ}
Let $(P_{n},\sigma)$ be a signed path. Then $r(P_{n},\sigma)=n-1$ if $n$ is odd and $r(P_{n},\sigma)=n$ if $n$ is even.
\end{lemma}
\noindent\begin{lemma}\label{le:2.2}\cite{FDD}
Let $(C_{n},\sigma)$ be a signed cycle, then
\begin{enumerate}[(a)]
\item if $(C_{n},\sigma)$ is balanced, then $r(C_{n},\sigma)=n-2$ if $n\equiv 0\;(mod\;4)$ and $r(C_{n},\sigma)=n$ otherwise.
\item if $(C_{n},\sigma)$ is unbalanced, then $r(C_{n},\sigma)=n-2$ if $n\equiv 2\;(mod\;4)$ and $r(C_{n},\sigma)=n$ otherwise.
\end{enumerate}
\end{lemma}

\noindent\begin{lemma}\label{le:2.3}\cite{HHL}
Let $(G,\sigma)$ be a signed graph, if $u$ is a pendant vertex of $G$ and $v$ is its unique neighbor, then $r(G,\sigma)=r((G,\sigma)-u-v)+2$.
\end{lemma}
\noindent\begin{lemma}\label{le:2.4}\cite{FDD}
Let $(G,\sigma)$ be a connected signed graph of order $n\geq 2$. Then
\begin{enumerate}[(a)]

\item $r(G,\sigma)=2$ if and only if $(G,\sigma)$ is a balanced complete bipartite graph.

\item $r(G,\sigma)=3$ if and only if $(G,\sigma)$ is a complete tripartite graph with partition $(V_{1}, V_{2}, V_{3})$, which satisfies that for each $i=1,2,3$, there exists a vertex $u_{i}\in V_{i}$ such that for every other vertex $v\in V_{i}$, either $N_{+}(v)=N_{+}(u_{i}), N_{-}(v)=N_{-}(u_{i})$ or $N_{+}(v)=N_{-}(u_{i}), N_{-}(v)=N_{+}(u_{i})$.
\end{enumerate}
\end{lemma}

 We use dotted line and solid line to denote negative edge and positive edge of a signed graph in  Figures 2, 3 and 4, respectively.

\noindent\begin{lemma}\label{le:2.5}\cite{ZY}
Let $(G,\sigma)$ be a signed bipartite graph. Then $r(G,\sigma)=4$ if and only if the reduced
graph $(H,\sigma)$ of $(G,\sigma)$ is switching equivalent to one of the signed graphs in Figure 2.
\end{lemma}

\noindent\begin{lemma}\label{le:2.6}\cite{ZY}
Let $(G,\sigma)$ be a signed nonbipartite graph and $(H,\sigma)$ be the reduced graph of $(G,\sigma)$. Then $r(G,\sigma)=4$ if and only if $(H,\sigma)$ is switching equivalent to one of the signed graphs in Figures 3 and 4.
\end{lemma}

\noindent\begin{lemma}\label{le:2.8}\cite{HHL}
Let $(G,\sigma)$ be a signed graph. If $(H,\sigma)$ is an induced subgraph of $(G,\sigma)$, then $r(G,\sigma)\geq r(H,\sigma)$.
\end{lemma}

\noindent\begin{lemma}\label{le:2.9}\cite{LW}
If $(G,\sigma)$ is a signed graph and different from a cycle, then $\eta(G,\sigma)\leq p(G)+2c(G)-1$, where $p(G)$ is the number of pendant vertices of $G$.
\end{lemma}



For a simple graph, we have the following lemma. Denote by $\bar{G}$ the complement graph of $G$.

\noindent\begin{lemma}\label{le:2.7}\cite{ZWT}
Let $G$ be a graph with girth $g$ and let $C$ be a shortest cycle in $G$. If there exists $u\in V(\bar{C})$ such that $|N_{C}(u)|\geq 2$, then $g=3$ or $4$.
\end{lemma}

\section{Bounds for the rank of $(G,\sigma)$ in terms of girth}
\quad In this section, we will  obtain some bounds for the rank of a signed graph in terms of its girth.
\noindent\begin{lemma}\label{le:3.1}
Let $(G,\sigma)$ be a connected signed graph and $(H,\sigma)$ be an induced subgraph of $(G,\sigma)$. If $r(G,\sigma)=r(H,\sigma)$ or $r(H,\sigma)+1$, then every vertex of $V(\bar{H})$ has a neighbor in $V(H)$.
\end{lemma}
\noindent\textbf{Proof.}
Suppose on the contrary that there are two vertices $u_{1}$ and $u_{2}$ in $G$ such that $d(u_{1},H)=1$ and $d(u_{2},H)=2$, where $u_{1}$ is adjacent $u_{2}$, respectively. By Lemma \ref{le:2.3}, one has that $$r((H,\sigma)+u_{1}+u_{2})=r(H,\sigma)+2\geq r(G,\sigma)-1+2=r(G,\sigma)+1,$$ since $(H,\sigma)+u_{1}+u_{2}$ is an induced subgraph of $(G,\sigma)$, we obtain a contradiction to Lemma \ref{le:2.8}.

\noindent\begin{theorem}\label{th:3.2}
Let $(G,\sigma)$ be a connected signed graph with girth $g$. Then
$r(G,\sigma)\geq g-2$, where the equality holds if and only if $(G,\sigma)$ is a signed graph of one of the followings.
\begin{enumerate}[(a)]
\item a balanced complete bipartite graph with $g=4$;
\item a balanced cycle with $g\equiv 0\;(mod\; 4)$;
\item an unbalanced cycle with $g\equiv 2\;(mod\; 4)$.
\end{enumerate}
\end{theorem}
\noindent\textbf{Proof.}
Let $(C,\sigma)$ be a shortest cycle in $(G,\sigma)$. Then by Lemmas \ref{le:2.2} and \ref{le:2.8}, we have that $r(G,\sigma)\geq r(C,\sigma)\geq g-2$. Now we prove the equality case.

\textbf{Sufficiency:} If $(G,\sigma)$ is a balanced complete bipartite graph with $g=4$, then by Lemma \ref{le:2.4}(a), we have that $$r(G,\sigma)=2=g-2.$$ If $(G,\sigma)$ is a balanced cycle with $g\equiv 0\;(mod\; 4)$ or an unbalanced cycle with $g\equiv 2\;(mod\; 4)$, then by Lemma \ref{le:2.2}, $$r(G,\sigma)=g-2.$$

\textbf{Necessity:} Because $r(G,\sigma)=g-2$, one has that $(C,\sigma)$ is a balanced cycle with $g\equiv 0\;(mod\; 4)$ or an unbalanced cycle with $g\equiv 2\;(mod\; 4)$; otherwise, by Lemmas \ref{le:2.2} and \ref{le:2.8}, we have that $$r(G,\sigma)\geq r(C,\sigma)=g,$$ a contradiction.
Now, we discuss two cases.

\textbf{Case 1.} Suppose $(C,\sigma)$ is a balanced cycle with $g\equiv 0\;(mod\; 4)$.

 \textbf{Subcase 1.1.}
 When $(G,\sigma)$ is a cycle, (b) of this theorem holds.

  \textbf{Subcase 1.2.}
 If $(G,\sigma)$ is not a cycle. Then by Lemma \ref{le:2.2}(a), $$r(G,\sigma)=g-2=r(C,\sigma).$$ So by Lemma \ref{le:3.1}, any vertex of $V(\bar{C})$ has a neighbor in $V(C)$. If for every vertex $u$ of $V(\bar{C})$, we have $|N_{C}(u)|\geq 2$, then by Lemma \ref{le:2.7}, $g=4$. So $$r(G,\sigma)=g-2=2,$$ by Lemma \ref{le:2.4}(a), $(G,\sigma)$ is a balanced complete bipartite graph. If there exists a vertex $u$ of $V(\bar{C})$ such that $|N_{C}(u)|= 1$, then by Lemmas \ref{le:2.1} and \ref{le:2.3}, $$r((C,\sigma)+u)=r(P_{g-1},\sigma)+2=g>g-2=r(G,\sigma),$$ a contradiction to Lemma \ref{le:2.8}. Combining with above, we can obtain that $(G,\sigma)$ is a balanced complete bipartite graph with $g=4$ or a balanced cycle with $g\equiv 0\;(mod\; 4)$.

\textbf{Case 2.} Suppose $(C,\sigma)$ is an unbalanced cycle with $g\equiv 2\;(mod\; 4)$, we have $g\geq 6$. Similarly to above, any vertex $u$ of $V(\bar{C})$ has a neighbor in $C$, by Lemma \ref{le:2.7}, $|N_{C}(u)|= 1$. Using the same method as in Case 1, we obtain $(G,\sigma)$ is an unbalanced cycle with $g\equiv 2\;(mod\; 4)$.

This completes the proof.  \quad $\square$

\noindent\begin{lemma}\label{le:3.3}
Let $(G,\sigma)$ be a connected signed graph with girth $g$. If $(G,\sigma)$ is not a balanced complete bipartite graph with $g=4$,
a balanced cycle with $g\equiv 0\;(mod\; 4)$ or an unbalanced cycle with $g\equiv 2\;(mod\; 4)$, then $r(G,\sigma)\geq g$.
\end{lemma}
\noindent\textbf{Proof.} By Theorem \ref{th:3.2}, we have that $$r(G,\sigma)\geq g-1.$$ If $r(G,\sigma)=g-1$, let $(C,\sigma)$ be a shortest cycle in $(G,\sigma)$, assume that $(C,\sigma)$ is neither an unbalanced cycle with $g\equiv 2\;(mod\; 4)$ nor a balanced cycle with $g\equiv 0\;(mod\; 4)$, then by Lemmas \ref{le:2.2} and \ref{le:2.8}, $$r(G,\sigma)\geq r(C,\sigma)=g>g-1=r(G,\sigma),$$ a contradiction. Now, we discuss two cases.

\textbf{Case 1.} Suppose $(C,\sigma)$ is an unbalanced cycle with $g\equiv 2\;(mod\; 4)$, so by Lemma \ref{le:2.2}(b), $$r(G,\sigma)=g-1=r(C,\sigma)+1,$$ then by Lemma \ref{le:3.1}, every vertex $u\in V(\bar{C})$ has a neighbor in $C$. Since $g\geq 6$, by Lemma \ref{le:2.7}, we have that for any vertex $u$ of $V(\bar{C})$, $|N_{C}(u)|= 1$. For some vertex $u_{1}$ of $V(\bar{C})$, by Lemmas \ref{le:2.1}, \ref{le:2.3} and \ref{le:2.8}, one has that $$g-1=r(G,\sigma)\geq r((C,\sigma)+u_{1})=r(P_{g-1},\sigma)+2=g,$$ a contradiction. So $(G,\sigma)$ is an unbalanced cycle with $g\equiv 2\;(mod\; 4)$, a contradiction.

\textbf{Case 2.} Suppose $(C,\sigma)$ is a balanced cycle with $g\equiv 0\;(mod\; 4)$. If $g=4$, $r(G,\sigma)=3$, then by Lemma \ref{le:2.4}(b), $(G,\sigma)$ is a complete tripartite graph, so $g=3$, a contradiction. If $g\geq 8$, using the same method as in Case 1, we can obtain that
$(G,\sigma)$ is a balanced cycle with $g\equiv 0\;(mod\; 4)$, a contradiction.

Thus $r(G,\sigma)\neq g-1$, then $r(G,\sigma)\geq g$.

This complete the proof.   \quad $\square$

\noindent\begin{lemma}\label{le:3.4}
Let $(G,\sigma)$ be a connected signed graph with girth $g$. If $r(G,\sigma)=g$ and $g$ is odd, then one of the following results is established:
\begin{enumerate}[(a)]
\item $(G,\sigma)$ is an odd cycle;
\item $(G,\sigma)$ is a complete tripartite graph ($g=3$) with partition $(V_{1}, V_{2}, V_{3})$, which satisfies that for each $i=1,2,3$, there exists a vertex $u_{i}\in V_{i}$ such that for every other vertex $v\in V_{i}$, either $N_{+}(v)=N_{+}(u_{i}), N_{-}(v)=N_{-}(u_{i})$ or $N_{+}(v)=N_{-}(u_{i}), N_{-}(v)=N_{+}(u_{i})$.
\end{enumerate}
\end{lemma}
\noindent\textbf{Proof.}
Because $g$ is odd, $g\geq 3$. Assume that $g=3$, that is $r(G,\sigma)=3$, then by Lemma \ref{le:2.4}(b), we have that $(G,\sigma)$ is a  complete tripartite graph with partition $(V_{1}, V_{2}, V_{3})$, which satisfies that for each $i=1,2,3$, there exists a vertex $u_{i}\in V_{i}$ such that for every other vertex $v\in V_{i}$, either $N_{+}(v)=N_{+}(u_{i}), N_{-}(v)=N_{-}(u_{i})$ or $N_{+}(v)=N_{-}(u_{i}), N_{-}(v)=N_{+}(u_{i})$.

If $g\geq 5$, let $(C,\sigma)$ be a shortest cycle of $(G,\sigma)$, then by Lemma \ref{le:2.2}, $$r(C,\sigma)=g=r(G,\sigma).$$ So by Lemmas \ref{le:2.7} and \ref{le:3.1}, for every vertex $u\in V(\bar{C})$, we have that $|N_{C}(u)|=1$. Let $u_{1}$ be a vertex in $V(\bar{C})$,
then by Lemmas \ref{le:2.1} and \ref{le:2.3},$$r((C,\sigma)+u_{1})=r(P_{g-1},\sigma)+2=g+1> r(G,\sigma),$$ a contradiction to Lemma \ref{le:2.8}.
Thus $(G,\sigma)$ is an odd cycle.

This completes the proof.   \quad $\square$

\noindent\begin{lemma}\label{le:3.5}
Signed bicyclic graph $(\theta(2m,2k,2l),\sigma)$ is nonsingular, where $m,k,l$ are positive integers and at most one of $m,k,l$ is 1.
\end{lemma}
\noindent\textbf{Proof.}
We show the underlying graph of the signed bicyclic graph $(\theta(2m,2k,2l),\sigma)$ in Figure 5.
Let $x_{1}$ and $y_{1}$ be the left and right $3$-degree vertices, respectively. Denote by $A(G,\sigma)$ the adjacent matrix of  $(\theta(2m,2k,2l),\sigma)$. We will prove that $A(G,\sigma)X=\textbf{0}$ only has zero solution. Suppose on the contrary that there exists a nonzero vector $\beta$ which is a solution of $A(G,\sigma)X=\textbf{0}$. Let $u$ be the element of $\beta$ corresponding to vertex $x_{1}$. Since $A(G,\sigma)\beta=\textbf{0}$, for each vertex $x$ of $(\theta(2m,2k,2l),\sigma)$, the sum of multiplying the elements of $\beta$ corresponding to the vertices adjacent to $x$ by the signs of the edges which are between the vertices adjacent to $x$ and vertex $x$ is $0$. Then we can obtain that the elements of $\beta$ for the three vertices
adjacent to vertex $y_{1}$ are equal to $u$ or $-u$. Thus for the vertex $y_{1}$, we have that $(\pm u)+(\pm u)+(\pm u)=\pm u=0$ or $\pm3u=0$, so $u=0$.

In a similar way, we also show that the elements of $\beta$ for the remaining vertices of $(\theta(2m,2k,2l),\sigma)$ are also $0$. Thus $\beta$ is the zero vector, a contradiction.
This proves that the signed bicyclic graph $(\theta(2m,2k,2l),\sigma)$ is nonsingular, where $m,k,l$ are positive integers and at most one of $m,k,l$ is 1.

This completes the proof.   \quad $\square$

\begin{figure}[htbp]
\centering
\includegraphics[scale=0.7]{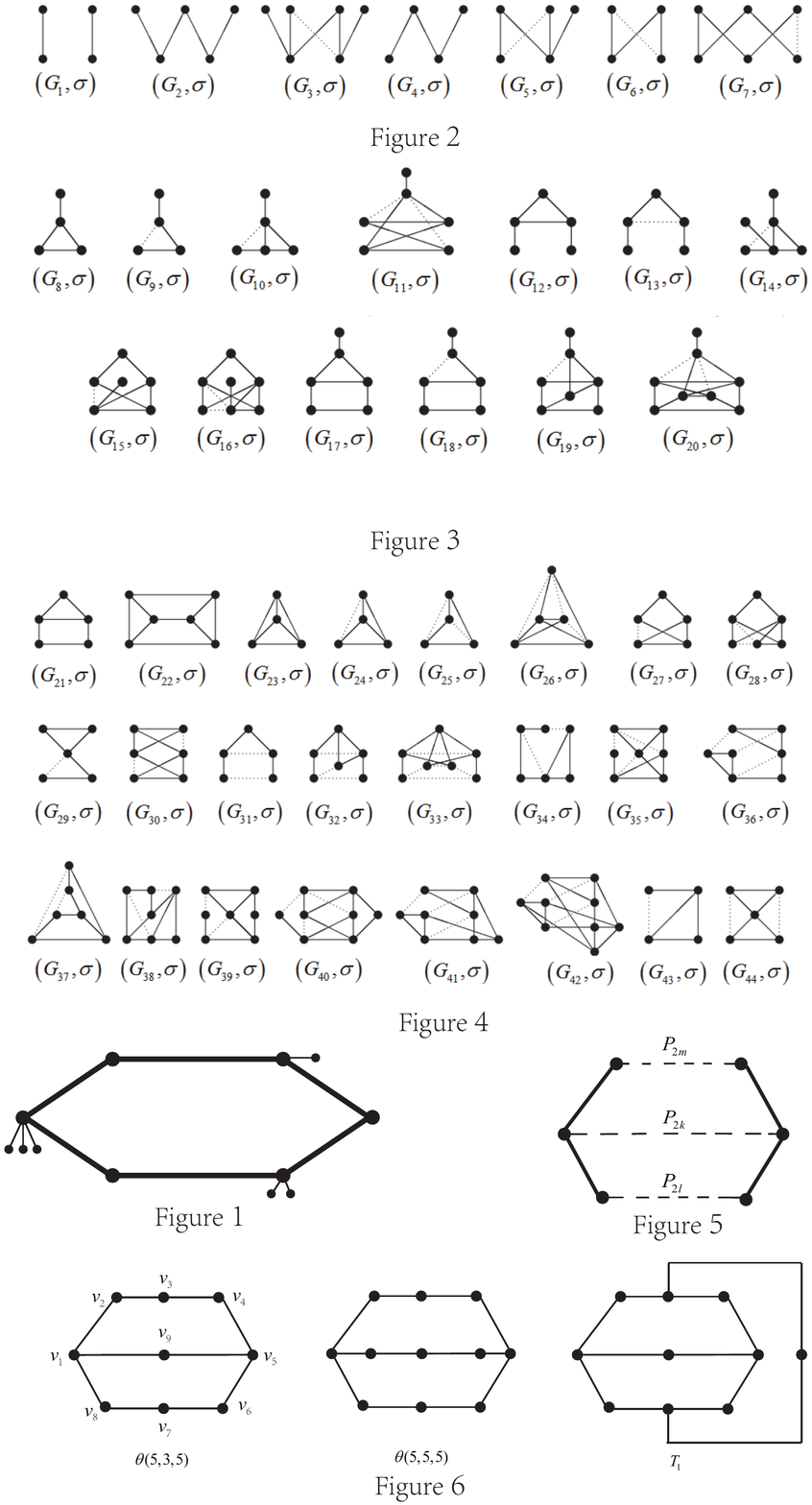}
\end{figure}

\noindent\begin{lemma}\label{le:3.6}
\begin{enumerate}[(a)]
\item Let $(G,\sigma)$ be the signed bicyclic graph. If  $(G,\sigma)$ is $(\theta(5,3,5),\sigma)$, then $r(G,\sigma)=g=6$ if and only if two cycles with order $6$ in $(G,\sigma)$ are both unbalanced.  If $(G,\sigma)$ is $(\theta(5,5,5),\sigma)$, then  $r(G,\sigma)=g=8$ if and only if $(G,\sigma)$ is a balanced signed bicyclic graph.

\item Let $(G,\sigma)=(T_{1},\sigma)$ be the signed tricyclic graph whose underlying graph is $T_{1}$ in Figure 6, then $r(G,\sigma)=g=6$ if and only if the cycles with order $6$ in $(G,\sigma)$ are all unbalanced.
\end{enumerate}
\end{lemma}

\noindent\textbf{Proof.}
 We use $v_{1}, v_{2},\ldots, v_{9}$ to label the vertices of $(\theta(5,3,5),\sigma)$ as Figure 6. Let $\beta$ be a solution of $A(G,\sigma)X=\textbf{0}$ and $w_{i}$ be the element of $\beta$ corresponding to the vertex $v_{i}$ ($i=1,2,\ldots,9$). Denote by $\alpha_{i}$ the row of $A(G,\sigma)$ corresponding to the vertex $v_{i}$ ($i=1,2,\ldots,9$). Since $A(G,\sigma)\beta=\textbf{0}$,  we have $$\beta=(w_{1},w_{2},w_{3},w_{4},w_{5},w_{6},w_{7},w_{8},w_{9})^{T}=(w_{1},w_{2},\pm w_{1},\pm w_{2},\pm w_{1},\pm w_{8},\pm w_{1},w_{8},\pm w_{2}\pm w_{8})^{T}.$$

\textbf{Necessity:} Let  $(G,\sigma)$ be $(\theta(5,3,5),\sigma)$. Since $r(G,\sigma)=g=6$, we have $\eta(G,\sigma)=3$, so $w_{1},w_{2}$ and $w_{8}$ are all nonzero, where $k_{1}w_{1}\neq k_{2}w_{2}\neq k_{3}w_{8}$ ($k_{1},k_{2}$ and $k_{3}$ are all nonzero). Then we have the cycle of order $8$ in $(G,\sigma)$ is balanced; otherwise, by Lemmas \ref{le:2.2} and \ref{le:2.8}, $$g=r(G,\sigma)\geq r(C_{8},\sigma)=g+2,$$ a contradiction.
Suppose on the contrary there exists a balanced cycle of order 6, without loss of generality, assume that it is the cycle $v_{1}v_{2}v_{3}v_{4}v_{5}v_{9}v_{1}$.   Now we consider two cases.

\textbf{Case 1.} The edges $v_{1}v_{2},v_{1}v_{9}$ have the same signs. Note that $\alpha_{i}\beta=0$ for each $i\in\{1,2,\ldots,9\}$.

\textbf{Subcase 1.1.} The edges $v_{2}v_{3},v_{3}v_{4}$ have the same signs. So by $\alpha_{3}\beta=0$, we have $w_{4}=-w_{2}$. Since the cycle $v_{1}v_{2}v_{3}v_{4}v_{5}v_{9}v_{1}$ is balanced, we have that the edges $v_{4}v_{5},v_{5}v_{9}$ have the same signs.  So by $\alpha_{1}\beta=0$, we have $w_{9}=-w_{2}\pm w_{8}$ and by $\alpha_{5}\beta=0$, we have $w_{9}=w_{2}\pm w_{8}$, then $w_{2}=0$ or $w_{2}=\pm w_{8}$, a contradiction.

\textbf{Subcase 1.2.} The edges $v_{2}v_{3},v_{3}v_{4}$ have the different signs. So by $\alpha_{3}\beta=0$, we have $w_{4}=w_{2}$. Since the cycle $v_{1}v_{2}v_{3}v_{4}v_{5}v_{9}v_{1}$ is balanced, we have that the edges $v_{4}v_{5},v_{5}v_{9}$ have the different signs. So by $\alpha_{1}\beta=0$, we have $w_{9}=-w_{2}\pm w_{8}$ and by $\alpha_{5}\beta=0$, we have $w_{9}=w_{2}\pm w_{8}$, then $w_{2}=0$ or $w_{2}=\pm w_{8}$, a contradiction.

\textbf{Case 2.} The edges $v_{1}v_{2},v_{1}v_{9}$ have the different signs. Using the same method as Case 1, we can obtain a contradiction.

Combining with above, one has that two cycles with order $6$ in $(G,\sigma)$ are both unbalanced.

\textbf{Sufficiency:} Since two cycles with order $6$ in $(G,\sigma)$ are both unbalanced, using the same method as above, one has that value of $w_{9}$ obtained by $\alpha_{1}\beta=0$ is the same as that obtained by $\alpha_{5}\beta=0$, that is $w_{9}=l_{1}w_{2}+l_{2}w_{8}\neq 0$ ($w_{2},w_{8},l_{1},l_{2}$ are all nonzero). Similarly, value of $w_{5}$ obtained by $\alpha_{9}\beta=0$ is the same as that obtained by $\alpha_{4}\beta=0$, and is also the same as that obtained by $\alpha_{6}\beta=0$, that is $w_{5}=l_{3}w_{1}\neq 0$ ($l_{3}\neq 0$). So we have that $$l_{1}w_{1}\neq l_{2}w_{2}\neq l_{3}w_{8}\neq 0.$$ So $\eta(G,\sigma)=3$, then $r(G,\sigma)=g=6$.

Using the same method as above, we can obtain that $r(\theta(5,5,5),\sigma)=g$ if and only if $(\theta(5,5,5),\sigma)$ is a balanced signed graph.

Using the same method as (a), we can obtain the (b) of this lemma.

This completes the proof.   \quad $\square$

Denote by $F_{1}$ the set of graphs that any graph in $F_{1}$ is obtained by identifying a vertex of the cycle $C_{n}$ with the center vertex $v$ of the star $S_{q}$(i.e., $V(C_{n})\cap V(S_{q})=v$). Let  $v_{1},v_{2},\ldots,v_{i}$ be the center vertex of stars $S_{q_{1}},S_{q_{2}},\ldots,S_{q_{i}}$, respectively. We use $F_{2}$ to denote the set of graphs that any graph in $F_{2}$ is obtained by identifying the vertices $u_{1},u_{2},\ldots,u_{i}$ ($2\leq i\leq n$) of the cycle $C_{n}$ with the center vertices $v_{1},v_{2},\ldots,v_{i}$, respectively (i.e., $V(C_{n})\cap V(S_{q_{j}})=v_{j}=u_{j}$, $j=1,2\ldots,i$).
\noindent\begin{lemma}\label{le:3.7}
Let $(G,\sigma)$ be a signed canonical unicyclic graph with girth $g$ and different from a
cycle. Then $r(G,\sigma)=g$ if and only if the following conditions hold:
\begin{enumerate}[(a)]
\item $g$ is even;
\item $G\in F_{1}$ or $G\in F_{2}$ and there is an odd number of vertices in the cycle of $G$ between any two consecutive star center vertices.
\end{enumerate}
\end{lemma}
\noindent\textbf{Proof.}
\textbf{Necessity:} Let $S_{m_{1}},S_{m_{2}},\ldots,S_{m_{k}}$ be $k$ ($k\geq1$) pendant stars and $(C,\sigma)$ be the unique cycle of $(G,\sigma)$. We assume that $v_{m_{1}},v_{m_{2}},\ldots,v_{m_{k}}$ are star center vertices of $S_{m_{1}},S_{m_{2}},\ldots,S_{m_{k}}$, respectively. If $g$ is odd, then by Lemmas \ref{le:2.1} and \ref{le:2.3}, $$r((C,\sigma)-v_{m_{1}}+V(S_{m_{1}}))=2+r(P_{g-1},\sigma)=g+1> g=r(G,\sigma),$$ a contradiction to Lemma \ref{le:2.8}. So $g$ is even. Suppose on the contrary that there exists two consecutive star center vertices such that there are an even number of vertices in $C$ between them, without loss of generality, let the two vertices be $v_{m_{1}}$ and $v_{m_{2}}$. Assume that there are $l_{1}$ ($l_{1}$ is even) vertices between $v_{m_{1}}$ and $v_{m_{2}}$, then by Lemmas \ref{le:2.1} and \ref{le:2.3}, $$r((C,\sigma)-v_{m_{1}}+V(S_{m_{1}})-v_{m_{2}}+V(S_{m_{2}}))=4+r(P_{g-2-l_{1}},\sigma)+r(P_{l_{1}},\sigma)=g+2> g=r(G,\sigma),$$ a contradiction to Lemma \ref{le:2.8}.

\textbf{Sufficiency:} Let $S_{m_{1}},S_{m_{2}},\ldots,S_{m_{k}}$ be $k$ ($k\geq1$) pendant stars and $(C,\sigma)$ be the unique cycle of $(G,\sigma)$. Let $n_{1},n_{2},\ldots,n_{k}$ be the orders of signed paths obtained by deleting the $k$ pendant stars from $(G,\sigma)$. Then by Lemmas \ref{le:2.1} and \ref{le:2.3}, we have that $$r(G,\sigma)=2k+\sum_{i=1}^{k} r(p_{n_{i}},\sigma)=2k+\sum_{i=1}^{k} (n_{i}-1)=k+\sum_{i=1}^{k} n_{i}=g.$$
This completes the proof.   \quad $\square$

\noindent\begin{lemma}\label{le:3.8}
Let $(G,\sigma)$ be a connected signed graph with girth $g\equiv 2\;(mod\;4)$. If $r(G,\sigma)=g$, then one of the following results is established:
\begin{enumerate}[(a)]
\item $(G,\sigma)$ is a balanced cycle with $g\equiv 2\;(mod\;4)$;
\item $(G,\sigma)$ is a signed canonical unicyclic graph;
\item $(G,\sigma)$ is a signed graph obtained by joining a vertex of an unbalanced cycle $(C_{g},\sigma)$ with $g\equiv 2\;(mod\;4)$ to the center of a signed star graph $(S_{k+1},\sigma)$ ($k\geq 1$);
\item $(G,\sigma)$ is $(\theta(5,3,5),\sigma)$ whose two cycles of order 6 are both unbalanced ($g=6$);
\item $G$ is tricyclic graph $T_{1}$ in Figure 6 and the cycles with order 6 of $(G,\sigma)$ are all unbalanced ($g=6$).
\end{enumerate}
\end{lemma}
\noindent\textbf{Proof.} If $(G,\sigma)$ is a cycle, then by Lemma \ref{le:2.2}, we have $(G,\sigma)$ is a balanced cycle with girth $g\equiv 2\;(mod\;4)$. We can obtain the (a) of this lemma.

Now we assume $(G,\sigma)$ is different from a cycle.
Let $(C,\sigma)$ be a shortest cycle of $(G,\sigma)$. Now we consider two cases.

\textbf{Case 1.} $(C,\sigma)$ is an unbalanced cycle. Because $g\geq 6$, then by Lemma \ref{le:2.7}, for every vertex $u\in V(\bar{C})$, we have that $|N_{C}(u)|\leq 1$. Denote by $U_{r}$ the set of all vertices $u$ in $V(\bar{C})$ such that $d(u,C)=r$. We can obtain that $U_{3}=\emptyset$; otherwise, choose a vertex $u_{3}\in U_{3}$, and let $u_{3}u_{2}u_{1}u_{0}$ be a shortest path in $G$ from $u_{3}$ to $C$, where $u_{2}\in U_{2}$, $u_{1}\in U_{1}$ and $u_{0}\in V(C)$. By Lemmas \ref{le:2.1} and \ref{le:2.3}, one has that $$r((C,\sigma)+u_{1}+u_{2}+u_{3})=4+r(P_{g-1},\sigma)=g+2> g=r(G,\sigma),$$ a contradiction to Lemma \ref{le:2.8}.

Observe that any two different vertices of $U_{1}$ are not adjacent; otherwise, there exists two vertices $v,v_{1}\in U_{1}$ such that $vv_{1}\in E(G)$. We assume that the neighbor vertices of $v,v_{1}$ in $C$ are $w,w_{1}$, respectively. Now adding the shorter of two paths in $C$ from $w$ to $w_{1}$ to the path $wvv_{1}w_{1}$, then we obtain a cycle of length at most $\frac{g}{2}+3$. Since $g$ is girth, one has that $\frac{g}{2}+3\geq g$, then $g\leq 6$. Because $g\geq 6$, then $g=6$. By Lemmas \ref{le:2.8} and \ref{le:3.5}, $$r(G,\sigma)\geq r((C,\sigma)+v+v_{1})=r(\theta(4,4,4),\sigma)=g+2>g,$$ a contradiction. If $U_{2}=\emptyset$, then $(G,\sigma)$ is a signed canonical unicyclic graph. We can obtain the (b) of this lemma.

Now we consider the following two subcases.

\textbf{Subcase 1.1.} $|U_{2}|\geq 2$. We have that for any two different vertices $m,l\in U_{2}$, $N_{U_{1}}(m)=N_{U_{1}}(l)$; otherwise, there are two distinct vertices $m_{2},l_{2}\in U_{2}$ and two distinct vertices $m_{1},l_{1}\in U_{1}$ such that $m_{2}m_{1}\in E(G)$ and $l_{2}l_{1}\in E(G)$. Assume $m_{2}l_{2}\notin E(G)$, then $m_{2}l_{1}\in E(G)$ and $m_{1}l_{2}\in E(G)$ do not both hold, otherwise $(G,\sigma)$ has a cycle $m_{2}l_{1}l_{2}m_{1}m_{2}$ with order 4, a contradiction to $g\geq 6$. So one or none of $m_{2}l_{1}\in E(G)$ and $m_{1}l_{2}\in E(G)$ holds. Then by Lemmas \ref{le:2.2}(b) and \ref{le:2.3}, $$r((C,\sigma)+m_{2}+l_{2}+m_{1}+l_{1})=r(C,\sigma)+4=g+2> g=r(G,\sigma),$$ a contradiction to Lemma \ref{le:2.8}. Thus $m_{2}l_{2}\in E(G)$, then we have that $m_{2}l_{1}\notin E(G)$ and $l_{2}m_{1}\notin E(G)$, otherwise $(G,\sigma)$ has a cycle of order 3, which contradicts to $g\geq 6$. Let the neighbor vertices of $m_{1}$ and $l_{1}$ in $C$ are $m_{0}$ and $l_{0}$, respectively. Since $(G,\sigma)$ contains no cycle of order 5, we have $m_{0}\neq l_{0}$. Now adding the shorter of two paths in $C$ from $m_{0}$ to $l_{0}$ to the path $l_{0}l_{1}l_{2}m_{2}m_{1}m_{0}$, then we obtain a cycle of length at most $\frac{g}{2}+5$. So one has that $\frac{g}{2}+5\geq g$, then $g\leq 10$. Because $g\equiv 2\; (mod\;4)$, then $g=6$ or $10$. By Lemma \ref{le:2.9}, signed bicyclic graphs without pendant vertices have nullity at most 3, so $\eta((C,\sigma)+m_{2}+l_{2}+m_{1}+l_{1})\leq 3$. Then by Lemma \ref{le:2.8}, $$g=r(G,\sigma)>r((C,\sigma)+m_{2}+l_{2}+m_{1}+l_{1})\geq g+4-3=g+1,$$ which is a contradiction.

Combining with above, we have $N_{U_{1}}(m)=N_{U_{1}}(l)$ for any two different vertices $m,l\in U_{2}$. So there has a vertex $s_{1}\in U_{1}$ such that $ss_{1}\in E(G)$ for all vertices $s\in U_{2}$.

Now we prove that $|U_{1}|=1$. Suppose on the contrary that there exists a vertex $o_{1}\in U_{1}$ with $o_{1}\neq s_{1}$. Let $s_{2},o_{2}\in U_{2}$ be two distinct vertices. Since $(G,\sigma)$ has no a cycle of order 4 and $N_{U_{1}}(s_{2})=N_{U_{1}}(o_{2})$, we have $s_{2}o_{1}\notin E(G)$ and $o_{2}o_{1}\notin E(G)$. Then by Lemmas \ref{le:2.1} and \ref{le:2.3}, $$r((C,\sigma)+s_{2}+s_{1}+o_{1})=4+r(P_{g-1},\sigma)=g+2> g=r(G,\sigma),$$ a contradiction to Lemma \ref{le:2.8}.
This shows that $|U_{1}|=1$. Thus $(G,\sigma)$ is a signed graph obtained by joining a vertex of an unbalanced cycle $(C_{g},\sigma)$ with $g\equiv 2\;(mod\;4)$ to the center of a signed star graph $(S_{k+1},\sigma)$ ($k\geq 2$). We obtain the (c) of this lemma.

\textbf{Subcase 1.2.} $|U_{2}|=1$, say $U_{2}=\{x\}$. If $|U_{1}|=1$, we also obtain the (c) of this lemma.

If $|U_{1}|\geq 2$, by the above proof, $U_{1}$ is an independent set. Let $y$ be the neighbor of $x$ in $U_{1}$ and $y_{1}\in U_{1}$ be different from $y$. If $xy_{1}\notin E(G)$, by Lemmas \ref{le:2.1} and \ref{le:2.3}, we have that $$r((C,\sigma)+x+y+y_{1})=r(P_{g-1},\sigma)+4=g+2> g=r(G,\sigma),$$ a contradiction to Lemma \ref{le:2.8}. Thus $xy_{1}\in E(G)$. Let the neighbor vertices of $y$ and $y_{1}$ in $(C,\sigma)$ be $z$ and $z_{1}$, respectively.

If $|U_{1}|=2$, by adding the shorter of two paths in $C$ from $z$ to $z_{1}$ to the path $z_{1}y_{1}xyz$, then we obtain a cycle of length at most $\frac{g}{2}+4$. So one has that $\frac{g}{2}+4\geq g$, then $g\leq 8$. Because $g\equiv 2\;(mod\;4)$, then $g=6$. Then $(G,\sigma)$ is $(\theta(5,3,5),\sigma)$, by Lemma \ref{le:3.6}, we have $r(G,\sigma)=g$ if and only if two cycles of order 6 in $(G,\sigma)$ are unbalanced. We obtain the (d) of this lemma.

If $|U_{1}|= 3$, suppose $y_{2}$ of $U_{1}$ is different from $y$ and $y_{1}$, where $xy\in E(G)$, using the same method as before, one has that $x y_{1}\in E(G)$ and $xy_{2}\in E(G)$. Let $z_{2}$ be the neighbor vertices of $y_{2}$ in $(C,\sigma)$. Then at least one of $d(z,z_{1}), d(z,z_{2})$ and $d(z_{1},z_{2})$ in $(C,\sigma)$ is at most $\lfloor\frac{g}{3}\rfloor$. Using the same method as before, we have that $\lfloor\frac{g}{3}\rfloor+4\geq g$, then $g\leq 6$. Because $g\equiv 2\;(mod\;4)$, then $g=6$. Thus by Lemma \ref{le:3.6}, $(G,\sigma)$ is a signed tricyclic graph whose underlying graph is $T_{1}$ in Figure 6 such that all cycles with order 6 in $(G,\sigma)$ are unbalanced. We can obtain the (e) of this lemma.

If $|U_{1}|\geq 4$, using the same method as before, we have that $\lfloor\frac{g}{4}\rfloor+4\geq g$, then $g\leq 5$, which is a contradiction.

\textbf{Case 2.} $(C,\sigma)$ is a balanced cycle. We have $U_{2}=\emptyset$; otherwise there exists a vertex $u\in U_{2}$ and let the neighbor of $u$ in $U_{1}$ be $v$. By Lemmas \ref{le:2.2}(a) and \ref{le:2.3}, $$r((C,\sigma)+v+u)=r(C,\sigma)+2=g+2> g=r(G,\sigma),$$  which contradicts to Lemma \ref{le:2.8}. Using the same methods as Case 1, we have $U_{1}$ is an independent set. Because $g\geq 6$, then by Lemma \ref{le:2.7}, we have $|N_{C}(u)|=1$ for all vertices $u\in V(\bar{C})$. So $(G,\sigma)$ is signed canonical unicyclic graph. We also obtain the (b) of this lemma.

This completes the proof. \quad $\square$

\noindent\begin{lemma}\label{le:3.9}
Let $(G,\sigma)$ be a connected signed graph with girth $g\equiv 0\;(mod\;4)$. If $r(G,\sigma)=g$, then one of the following results is established:
\begin{enumerate}[(a)]
\item the reduced graph of $(G,\sigma)$ is switching equivalent one of signed graphs in Figure 2 except for $(G_{1},\sigma)$ ($g=4$);
\item $(G,\sigma)$ is a unbalanced cycle with girth $g\equiv 0\;(mod\;4)$;
\item $(G,\sigma)$ is a signed canonical unicyclic graph;
\item $(G,\sigma)$ is a signed graph obtained by joining a vertex of a balanced cycle $(C_{g},\sigma)$ with $g\equiv 0\;(mod\;4)$ to the center of a signed star graph $(S_{k+1},\sigma)$ ($k\geq 1$);
\item $(G,\sigma)$ is a balanced signed bicyclic graph $(\theta(5,5,5),\sigma)$ ($g=8$).
\end{enumerate}
\end{lemma}
\noindent\textbf{Proof.} If $g=4$, then $r(G,\sigma)=4$. By Lemmas \ref{le:2.5} and \ref{le:2.6}, the reduced graph of $(G,\sigma)$ is switching equivalent to one of the signed graphs in Figures 2, 3 and 4. Since $(G,\sigma)$ is connected and $g=4$, then the reduced graph of $(G,\sigma)$ is switching equivalent to one of the signed graphs in Figure 2 except for $(G_{1},\sigma)$.
If $g\geq 8$, using the same method as the proof of Lemma \ref{le:3.8}, one has that $(G,\sigma)$ is of either (b),(c),(d) or (e).

This completes the proof. \quad $\square$

\noindent\begin{theorem}\label{th:3.10}
Let $(G,\sigma)$ be a connected signed graph with girth $g$. Then $r(G,\sigma)=g$ if and only if $(G,\sigma)$ is a graph of one of the followings:
\begin{enumerate}[(a)]
\item a signed odd cycle;
\item a balanced cycle with $g\equiv 2\;(mod\;4)$ or an unbalanced cycle with $g\equiv 0\;(mod\;4)$;
\item a complete tripartite graph ($g=3$) with partition $(V_{1}, V_{2}, V_{3})$, which satisfies that for each $i=1,2,3$, there exists a vertex $u_{i}\in V_{i}$ such that for every other vertex $v\in V_{i}$, either $N_{+}(v)=N_{+}(u_{i}), N_{-}(v)=N_{-}(u_{i})$ or $N_{+}(v)=N_{-}(u_{i}), N_{-}(v)=N_{+}(u_{i})$;
\item a signed canonical unicyclic graph with an even cycle such that either $G\in F_{1}$ or $G\in F_{2}$ and there is an odd number of vertices in the cycle of $G$ between any two consecutive star center vertices;
\item a signed graph obtained by joining the center of a signed star graph $(S_{k+1},\sigma)$ ($k\geq 1$) to either a vertex of an unbalanced cycle $(C_{g},\sigma)$ with $g\equiv 2\;(mod\;4)$ or a balanced cycle $(C_{g},\sigma)$ with $g\equiv 0\;(mod\;4)$;
\item a signed graph whose reduced graph is switching equivalent to one of Figure 2 except for $(G_{1},\sigma)$ ($g=4$);
\item a balanced signed bicyclic graph $(\theta(5,5,5),\sigma)$ ($g=8$) or $(\theta(5,3,5),\sigma)$ whose cycles with order 6 are both unbalanced ($g=6$);
\item a signed tricyclic graph whose underlying graph is $T_{1}$ in Figure 6 such that all cycles with order 6 of $(G,\sigma)$ are unbalanced ($g=6$).
\end{enumerate}
\end{theorem}
\noindent\textbf{Proof.}
\textbf{Necessity:} By Lemmas \ref{le:3.4}, \ref{le:3.7}, \ref{le:3.8} and \ref{le:3.9}, we have that the results hold.

\textbf{Sufficiency:} If $(G,\sigma)$ is of (a) or (b), by Lemma \ref{le:2.2}, one has that $r(G,\sigma)=g$. If $(G,\sigma)$ is of (c)(resp. (d)), by Lemma \ref{le:2.4}(b)(resp. \ref{le:3.7}), we have that $r(G,\sigma)=g$. If $(G,\sigma)$ is of (e), by Lemmas \ref{le:2.2} and \ref{le:2.3}, $$r(G,\sigma)=2+r(C_{g},\sigma)=g.$$ If $(G,\sigma)$ is of (f), the equality $r(G,\sigma)=g$ follows from Lemma \ref{le:2.5}. If $(G,\sigma)$ is of (g) or (h), by Lemma \ref{le:3.6}, we have that $r(G,\sigma)=g$.

This completes the proof. \quad $\square$



\end{document}